\documentclass[a4paper,11pt]{article}
\usepackage{graphicx}
\usepackage{amssymb}
\usepackage{latexsym, bm}
\usepackage{multicol}
\usepackage{enumerate}
\usepackage{indentfirst}
\usepackage{amssymb,amsfonts}
\usepackage{amsmath}
\usepackage{setspace}
\usepackage{CJK}

\usepackage{amssymb}

\newtheorem{theorem}{Theorem}
\newtheorem{lemma}{Lemma}
\newtheorem{fact}{Fact}
\newtheorem{corollary}{Corollary}

\newtheorem{problem}{Problem}

\newtheorem{claim}{Claim}

\title{ Cycle Ramsey numbers for random graphs\thanks{Supported in part by NSFC(11671088) and NSFFP(2016J01017).}
}
\author{Meng Liu$^1$,\,\,\,Yusheng Li$^2$,\,\,\,  Qizhong Lin$^3$\footnote{Corresponding Author.},\,\,\, and\,\,\,Chunlin You$^3$
\vspace*{0.3cm}\\
{\small $^1$ School of Mathematical Sciences, Anhui University, Hefei, Anhui 230601, China}\\
{\small$^2$ School of Mathematical Sciences, Tongji University, Shanghai 200092, China} \\
{\small $^3$ Center for Discrete Mathematics, Fuzhou University, Fuzhou, Fujian 350108, China} \\
{\small\em (liumeng@ahu.edu.cn (M. Liu),\, li\_yusheng@tongji.edu.cn (Y. Li),}
\\{\small\em linqizhong@fzu.edu.cn (Q. Lin), chunlin\_you@163.com (C. You))}
}

\date{}
\begin{spacing}{1.0}
\begin{document}
\maketitle
\begin{abstract}
Let $C_{n}$ be a cycle of length $n$.
As an application of Szemer\'{e}di's regularity lemma, {\L}uczak ($R(C_n,C_n,C_n)\leq (4+o(1))n$, J. Combin. Theory Ser. B, 75 (1999), 174--187) in fact established that $K_{(8+o(1))n}\to(C_{2n+1},C_{2n+1},C_{2n+1})$.
In this paper, we strengthen several results involving cycles.

Let $\mathcal{G}(n,p)$ be the random graph. We prove that for fixed $0<p\le1$, and integers $n_1$, $n_2$ and $n_3$ with $n_1\ge n_2\ge n_3$, it holds that for any sufficiently small $\delta>0$, there exists an integer $n_0$ such that for all integer $n_3>n_0$, we have a.a.s. that
\begin{align*}
\mathcal{G}((8+\delta)n_1,p) \to (C_{2n_1+1},C_{2n_2+1},C_{2n_3+1}).
\end{align*}
Moreover, we prove that for fixed $0<p\le1$ and integers $n_1\ge n_2\ge n_3>0$ with same order, i.e. $n_2=\Theta(n_1)$ and $n_3=\Theta(n_1)$, we have a.a.s. that
\begin{align*}
\mathcal{G}(2n_1+n_2+n_3+o(1)n_1,p) \to (C_{2n_1},C_{2n_2},C_{2n_3}).
\end{align*}
Similar results for the two color case are also obtained.

\medskip


\medskip


\end{abstract}

\section{Introduction}
For graphs $F$, $G_1,\dots,G_k$, we write $F \to (G_1,\dots,G_k)$ if every edge coloring of $F$ in $k$ colors
yields a monochromatic graph $G_i$ for some $1\le i\le k$. The Ramsey number $R(G_1,\dots,G_k)$ is the minimum integer $N$ such that
$
K_N \to (G_1,\dots,G_k).
$
 When $G_i=G$ for all $1\le i\le k$, we denote $F \to (G_1,\dots,G_k)$ by $F \to (G)_k$ and $R(G_1,\dots,G_k)$ by $R_k(G)$. We refer the reader to the book by Graham, Rothschild and Spencer \cite{graham-roth-spen} for an overview and a survey by Conlon, Fox and Sudakov \cite{cfs} for many recent developments.

The 2-color Ramsey numbers of cycles $R(C_{n_1},C_{n_2})$ are studied by several researchers, we refer the reader to Bondy and Erd\H{o}s \cite{B-e}, Faudree and Schelp \cite{Faudree-1974}, and independently, Rosta \cite{Rosta-1973} (see also \cite{k-r}).
A few years later Erd\H{o}s et al. \cite{er} found the value of $R(C_{n_1},C_{n_2},C_{n_3})$
and $R(C_{n_1},C_{n_2},C_{n_3},C_{n_4})$ when one of the cycles is much longer than the
others.
An elegant result by {\L}uczak \cite{Luczak-1999} proved that
\[
R_{3}(C_{2n+1} )=(8 + o(1))n.
\]
In that paper {\L}uczak introduced a technique that uses the regularity method to reduce problems about paths and cycles to problems about matchings that are contained in a connected component.
This technique has become fairly standard in the area, many results appeared, see e.g. \cite{Benevides-2009,Davies-2017,Figaj--2007,f-l,Gyarfas-2007,Gyarfas-M-2007,Kohayakawa-2005,Gsark-2016}.
In particular, Figaj and {\L}uczak \cite{Figaj--2007,f-l} obtained that for $\alpha_1\ge\alpha_2\ge\alpha_3>0$,
\begin{align*}
&R(C_{2\lfloor\alpha_1n\rfloor},C_{2\lfloor\alpha_2n\rfloor},C_{2\lfloor\alpha_3n\rfloor})=(2\alpha_1+\alpha_2+\alpha_3+o(1))n,
\\
&R(C_{2\lfloor\alpha_1n\rfloor+1},C_{2\lfloor\alpha_2n\rfloor+1},C_{2\lfloor\alpha_3n\rfloor+1})
=(8\alpha_1+o(1))n.
\end{align*}
Generally, a recent result by Jenssen and Skokan \cite{Jenssen-Skokan-2016}
proved that for every  $k$ and sufficiently large $n$,
\[
{R_k}({C_{2n+1}}) = {2^{k}}n + 1.
\]
This resolves a conjecture of Bondy and Erd\H{o}s \cite{B-e}.
In particular for $k=3$, the conjecture has been confirmed by Kohayakawa, Simonovits, and Skokan in an earlier paper \cite{Kohayakawa-2005}.
Also,  for every fixed integer $k\geq 4$ and sufficiently large even $n$,
 \[
(2k-2 + o(1))n\le{R_k}({C_{2n}}) \le (2k-1 + o(1))n,
 \]
 where the lower bound due to Sun et al. \cite{sun} while the upper bound
by Knierim and Su \cite{Knierim-2018}.

However, very little known about the Ramsey numbers for random graphs.
The {\em random graph} $\mathcal{G}(n,p)$ is the random graph $G$ with vertex set $[n]:=\{1,2,\dots,n\}$
in which each pair $\{i,j\}\in {[n]\choose 2}$ appears independently as an edge in $G$ with probability $p$.
We say an event $E_n$ occurs asymptotically almost surely (a.a.s.) if $\mathop {\lim }\limits_{n \to \infty } \Pr ({E_n}) = 1$.
Let $P_n$ be a path with $n$ vertices.
Applying the sparse regularity lemma established by Kohayakawa and R\"{o}dl \cite{koha97,Ko-r}, Letzer \cite{letzer-2016} showed that a.a.s. $\mathcal{G}(n,p) \to {P_{(2/3 - o(1) )n}}$, provided that $pn \to \infty $, which is optimal.
By using the idea of Figaj and {\L}uczak \cite{Figaj--2007}, Dudek and Pra{\L}at \cite{dudek-2018} proved that a.a.s.
$\mathcal{G}(n,p) \to {({P_{(1/2 - o(1))n}})_3}$ as $pn \to \infty $. For multicolor $k\ge3$, Dellamonica et al. \cite{dell}  proved that a.a.s.
$\mathcal{G}(n,p) \to {({P_{(1/k - o(1) )n}})_k}$, provided that $pn \to \infty $.

In this paper, we consider cycle Ramsey numbers for random graphs.

\begin{theorem}\label{odd-c-3}
Let $0<p\le1$ be fixed, and $n_1$, $n_2$ and $n_3$ be integers with $n_1\ge n_2\ge n_3$. Then for any sufficiently small positive real $\delta$, there exists an integer $n_0$ such that for all integer $n_3>n_0$, we have a.a.s. that
\[
\mathcal{G}((8+\delta)n_1,p) \to (C_{2n_1+1},C_{2n_2+1},C_{2n_3+1}),\;\text{and}\;\;
\mathcal{G}((4+\delta)n_1,p) \to (C_{2n_1+1},C_{2n_2+1}).
\]
\end{theorem}

For any integers $n_1$, $n_2$ and $n_3$ with $n_1\ge n_2\ge n_3$, we can easily see the lower bound $R(C_{2n_1+1},C_{2n_2+1},C_{2n_3+1})>8n_1$.
This together with Theorem \ref{odd-c-3} yield the following result which is slightly stronger than that by {\L}uczak \cite[Theorem 1]{Luczak-1999}, and Figaj and {\L}uczak \cite[Theorem 1 (iv)]{f-l}.
\begin{corollary}
For all large integers $n_1$, $n_2$ and $n_3$ with $n_1\ge n_2\ge n_3$, we have
\[
R(C_{2n_1+1},C_{2n_2+1},C_{2n_3+1})=(8+o(1))n_1.
\]
\end{corollary}

For even cycle, we have a similar result as follows.
\begin{theorem}\label{c-3}
Let $0<p\le1$ be fixed, and let $n_1\ge n_2\ge n_3>0$ be integers with same order, we have a.a.s. that
\begin{align*}
\mathcal{G}(2n_1+n_2+n_3+o(1)n_1,p) &\to (C_{2n_1},C_{2n_2},C_{2n_3}),\;\;\text{and}
\\ \mathcal{G}(2n_1+n_2+o(1)n_1,p) &\to (C_{2n_1},C_{2n_2}).
\end{align*}
\end{theorem}

Note that for integers $n_1$, $n_2$ and $n_3$ with $n_1\ge n_2\ge n_3$,
we can see the lower bound $R(C_{2n_1},C_{2n_2},C_{2n_3})\ge 2n_1+n_2+n_3-2$.
Indeed, let $K_N$ be defined on $V=\cup_{i=1}^3V_i$, where $V_1, V_2$ and $V_3$ are disjoint with $|V_1| = 2n_1-1$, $|V_2| = n_2-1$, $|V_1| = n_3-1$. We color the edges inside $V_1$ with the first color, all the edges adjacent to vertices from $V_2$ with the second color, and all the remaining edges of $K_N$ with the third color. It is easy to check that there is no monochromatic $C_{2n_i}$ as desired.
This together with Theorem \ref{c-3} yield the following result by  Figaj and {\L}uczak \cite[Theorem 1]{Figaj--2007}.
\begin{corollary}
For $\alpha_1\ge\alpha_2\ge\alpha_3>0$, we have
\begin{align*}
R(C_{2\lfloor\alpha_1n\rfloor},C_{2\lfloor\alpha_2n\rfloor},C_{2\lfloor\alpha_3n\rfloor})=(2\alpha_1+\alpha_2+\alpha_3+o(1))n.
\end{align*}
\end{corollary}

\section{Regularity lemma}

Let $A$ be a set of positive integers and $A_n=A\cap \{1,...,n\}$. In the 1930s, Erd\H{o}s and Tur\'an conjectured that if $\overline{\lim}\frac{|A_n|}{n}>0$,
then $A$ contains arbitrarily long arithmetic progressions.
The conjecture for the arithmetic progressions of length 3 was proved by Roth \cite{roth53,roth54}.
The full conjecture was proved by Szemer\'edi \cite{sze75} with a deep and complicated combinatorial argument.
In the proof, Szemer\'edi used a result, which is now called the bipartite regularity lemma.
For the general regularity lemma, see Szemer\'edi \cite{Sze-1978}.
The lemma has became a totally new tool in extremal graph theory. For many applications, we refer the reader to the survey of Koml\'os and Simonovits \cite{Kom-Simonovits-2005} and related references.

Let $G=(V,E)$ be a graph. If $U,W\subseteq V(G)$ are nonempty disjoint sets,
we write $e_G(U,W)$ for the edges between $U$ and $W$ in $G$, and call

\[d_G(U,W)=\frac{e_G(U,W)}{|U||W|}
\]
the density of the pair $(U,W)$.

Let $\epsilon >0$, a pair $(U,W)$ of nonempty disjoint sets
$U,W\subseteq V(G)$ is called $\epsilon$-regular if

\[
|d_G(X,Y)-d_G(U,W)|\le \epsilon
\]
for every $X\subseteq U, Y\subseteq W$ such that $|X|\ge \epsilon |U|$ and $|Y|\ge \epsilon |W|$.
The pair $(U,W)$ is called $(\epsilon,d)$-regular if it is $\epsilon$-regular and $d_G(U,W)\geq d$. Moreover, $(U,W)$ is called $(\epsilon,d)$-super-regular if it is $\epsilon$-regular and $\deg_G(u,W)>d|W|$ for all $u\in U$ and $\deg_G(w,U)>d|U|$ for all $w\in W$, where $\deg_G(u,W)$ is the number of neighbors lie in $W$ of a vertex $u$ in $G$.

Let us list several useful properties for regularity pairs.

\begin{fact}\label{pre}
 Let $(U,W)$ be an $(\epsilon,d)$-regular pair, and $Y\subseteq W$ with $|Y|\ge \epsilon |W|$. Then there
exists a subset $U'\subseteq U$ with $|U'|\ge (1-\epsilon)|U|$,
each vertex in $U'$ is adjacent to at least $(d-\epsilon)|Y|$
vertices in $Y$.
\end{fact}

 \begin{fact}\label{hdt}
 Let $(U,W)$ be an $\epsilon$-regular pair in $G$, and let $X\subseteq U$ and $Y\subseteq W$ with $|X|\geq \gamma|U|$ and $|Y|\geq \gamma|W|$ for some $\gamma>\epsilon$. Then $(X,Y)$ is $\epsilon'$-regular such that $|d_G(U,W)-d_G(X,Y)|<\epsilon$, where $\epsilon'=\max\{\epsilon/\gamma, 2\epsilon\}$.
\end{fact}

The following property tells that any regular pair has a large subgraph which is super-regular, and here we include a proof.
\begin{fact}\label{supreg}
 For $0<\epsilon<1/2$ and $d\leq1$, if $(U,W)$ is $(\epsilon,d)$-regular with $|U|=|W|=m$ then there exist $U_1\subseteq U$ and $W_1\subseteq W$ with $|U_1|=|W_1|=(1-\epsilon)m$ such that $(U_1,W_1)$ is $(2\epsilon,d-2\epsilon)$-super-regular.
 \end{fact}
\noindent{\em Proof.}
Let $X\subseteq U$ consists of vertices with at most $(d-\epsilon)|W|$ neighbors in $W$. Since $e(X,W)\leq |W|(d-\epsilon)|W|$ we have  $|d(X,W)-d|\geq \epsilon$.
According to the definition of $\epsilon$-regular, we know that $|X|<\epsilon m$.
Similarly, let $Y\subseteq W$ consists of vertices with at most $(d-\epsilon)|U|$ neighbors in $U$, we know that $|Y|<\epsilon m$.

Take $U_1\subseteq U\setminus X$ and $W_1\subseteq W\setminus Y$ with $|U_1|=|W_1|=(1-\epsilon)m$.
Clearly, each vertex of $U_1$ has at least $(d-\epsilon)m-\epsilon m=(d-2\epsilon)m$ neighbors in $W_1$, and similarly
each vertex of $W_1$ has at least $(d-2\epsilon)m$ neighbors in $U_1$.
On the other hand, for any subset $S\subseteq U_1$ and $T\subseteq W_1$.
If $|S|\ge 2\epsilon|U_1|$ and $|T|\ge 2\epsilon|W_1|$, then clearly $|S|\ge \epsilon m$ and $|T|\ge \epsilon m$.
From the fact that $(U,W)$ is $(\epsilon,d)$-regular, we have
\[
|d(S,T)-d(U_1,W_1)|\le|d(S,T)-d(U,W)|+|d(U_1,W_1)-d(U,W)| <2\epsilon.
\]
This completes the proof.\hfill$\Box$

In this note,  we use the following version of regularity lemma \cite{Sze-1978}.

\begin{lemma}\label{reg}(Regularity lemma)
For $\epsilon >0$ and integer $t_0\ge 1$, there exists $T_0=T_0(\epsilon,t_0)$ such that the following holds.
For all graphs $G_1$, $G_2$ and $G_3$ with the same vertex set $V$ and $|V|\ge t_0$, there exists a partition
$\{V_0,V_1,\dots,V_k\}$ of $V(G)$ such that $t_0\le t\le T_0$ and

\medskip
(1) $|V_0|<\epsilon n, |V_1|=|V_2|=\dots=|V_t|$;

(2) all but at most $\epsilon t^2$ pairs $(V_i,V_j)$, $i,j\in [t]$, are $\epsilon$-regular for $G_1$, $G_2$ and $G_3$.
\end{lemma}

\section{Proofs of main results}

In this section, we will give proofs for our results by using regularity method together with probability method.
Let us list several lemmas we will use to establish our main results.
We first recall a result by {\L}uczak \cite[Lemma 9]{Luczak-1999}.

\begin{lemma}({\L}uczak \cite{Luczak-1999})\label{Luczak-1999}
For every $0<\eta< 10^{-5}$ and $t\ge\exp(\eta^{-50})$ the following holds.
If $H$ is a graph with $T\ge4(1+\eta)t$ vertices and at least $(1-\eta^3){T\choose 2}$ edges,
then each 3-coloring of edges of $G$ leads to a monochromatic odd cycle of length at least
$(1+\eta/10)t$.
\end{lemma}

The next lemma by Benevides and Skokan \cite[Lemma 10]{Benevides-2009}, independently by Figaj and {\L}uczak \cite[Lemma 5]{Figaj--2007} about regularity pairs is a slightly stronger version of  \cite[Claim 3]{Luczak-1999} by {\L}uczak. For the proof, one can easily follow that of \cite{Luczak-1999}.

\begin{lemma}(Benevides and Skokan \cite{Benevides-2009}, Figaj and {\L}uczak \cite{Figaj--2007})\label{ben-sk}
For every $0<\beta< 1$, there exists an $n_0$ such that for every $n > n_0$ the following holds: Let $G$ be a bipartite graph with bipartition $V(G)=V_1\cup V_2$ such that $|V_1|=|V_2|=n$. Furthermore, let the pair $(V_1,V_2)$ be $\epsilon$-regular with density at least $\beta/4$ for some $\epsilon$ satisfying $0<\epsilon<\beta/100$. Then for every $\ell, 1\le \ell \le n-5\epsilon n/\beta$, and for every pair of vertices $v'\in V_1, v''\in V_2$ satisfying $\deg(v'),\deg(v'')\ge \beta n/5$, $G$ contains a path of length $2\ell+1$ connecting $v'$ and $v''$.
\end{lemma}

\noindent
{\em 3.1 \; Proof of Theorem \ref{odd-c-3}}

\bigskip

We first prove the three color case. Let $0<p\le1$ be fixed, and let
\begin{align}\label{eta-ep}
&0<\eta< 10^{-5},\;\;0<\epsilon\le\min\left\{\eta^3,{p^2}/{4}\right\},\;\;\text{and}\;\;t_0>4(1+\eta)\eta^{-50}.
\end{align}
Furthermore, denote $n=n_1$ for convenience and let
\begin{align}\label{alp}
\alpha=10\eta,\;\;\text{and}\;\;N=(8+\alpha)n.
\end{align}
We will show that a.a.s. for every 3-edge coloring of $G=\mathcal{G}(N,p)$ with vertex set $V$, there is a monochromatic odd cycle of length $2n_\lambda+1$ on the $\lambda$th color for $\lambda=1,2,3$.

Consider a 3-edge coloring $(G_1, G_2, G_3)$ of $G$.
By Lemma \ref{reg}, there exists $T_0=T_0(\epsilon,t_0)$ such that there is a partition $\{V_0,V_1,\dots,V_t\}$ of $V(G)$ satisfying $t_0\le t\le T_0$ and
\medskip

(1) $|V_0|<\epsilon N, |V_1|=|V_2|=\dots=|V_t|$;

(2) all but at most $\epsilon t^2$ pairs $(V_i,V_j)$, $i,j\in [t]$, are $\epsilon$-regular for $G_1$, $G_2$ and $G_3$.

\begin{claim}\label{claim-1}
For every two disjoint subsets $U$, $W\subseteq V$ of sizes much larger than $\ln N$, we have a.a.s.
\[
d(U,W)\geq p/2.
\]
\end{claim}
{\em Proof.} Let $|U|=k_1$ and $|W|=k_2$. By Chernoff's Inequality (see e.g. \cite{Alon-1992,Chernoff-1952}),
we have for some constant $c>0$ that $\Pr (d(U,W)<p/2 ) < \exp(- pck_1k_2)$ as the expectation of the density of $d(U,W)$ is $p$.
Let $X$ be the number of such pairs $(U,W)$ satisfying $d(U,W)<p/2.$
Then the expectation of $X$ satisfies
$$E(X)\leq {N\choose k_1}{N-k_1 \choose k_2}\cdot \exp\left(-pc{k_1k_2}\right),$$
which will tend to zero as $k_1$ and $k_2$ are much larger that $\ln N$ and $p$ fixed.
Therefore, by Markov's Inequality, it follows that $\Pr(X\geq 1)<E(X)\to 0$ as $N\to \infty$ as desired.\hfill $\Box$

By Claim \ref{claim-1}, we suppose that $d_{G}(V_i,V_j)\geq p/2$ for every  $1\leq i \leq j \leq t$.
Let $H$ be the auxiliary graph on vertex set $[s]$,
where $ij$ is an edge of $H$ if and only if  $(V_i,V_j)$ is $\epsilon$-regular.
Since the partition $V_0,V_1,\dots,V_t$ is $\epsilon$-regular, the  number of edges in $H$
is at least $(1-\epsilon){t\choose 2}.$
We color an edge $ij$ in $H$ red if $d_{G_1}(V_i,V_j)\geq p/6$, blue if $d_{G_2}(V_i,V_j)\geq p/6$ and green if $d_{G_3}(V_i,V_j)\geq p/6$.

From Lemma \ref{Luczak-1999}, we have that there is a monochromatic odd cycle $C$ of length $s$ which is at least
$(1+\eta/10)t/4(1+\eta)$.
Without loss of generality, suppose that the odd cycle $C$ with vertex set $[s]$ in $G_\lambda$ for some $\lambda=1,2,3$, i.e. for $i=1,2,\dots,s$,

\medskip
(a) $(V_i,V_{i+1})$ is $\epsilon$-regular, and

\smallskip
(b) $d_{G_\lambda}(V_i,V_{i+1})\geq p/6$, where the indices are taken modulo $s$.

\medskip
In the following, we aim to show that $G_\lambda$ contains an odd cycle of length $2n_\lambda+1$.
From Fact \ref{supreg}, for odd $i=1,3,\dots,s-2$, we can take $V_i'\subseteq V_i$
with $|V_i'|\ge (1-\epsilon)|V_i|$ and $|V_{i+1}'|\ge (1-\epsilon)|V_{i+1}|$ such that $(V'_i,V'_{i+1})$ is $(2\epsilon,p/6-2\epsilon)$-super-regular.
Let
\begin{align*}
m=(1-\epsilon)^2\frac{N}{t},\;\;\epsilon'=2\epsilon,\;\;\text{and}\;\;\beta=2p/3-8\epsilon.
\end{align*}
Then we have that for every pair of vertices $v_i'\in V_i'$, $v_{i+1}'\in V_{i+1}'$ satisfying
\[
\deg_{G_\lambda}(v_i',V_{i+1}')\ge (p/6-2\epsilon)|V_{i+1}'|\ge\beta m/5,\;\;
\text{and}\;\;\deg_{G_\lambda}(v_{i+1}',V_i')\ge \beta m/5.
\]
Therefore, by Lemma \ref{ben-sk}, for every $\ell$, $1\le \ell \le m-5\epsilon m/\beta$,
and for every pair of vertices $v_i'\in V_i'$ and $v_{i+1}'\in V_{i+1}'$, $G_\lambda$ contains a path of length $2\ell+1$ connecting $v_i'$ and $v_{i+1}'$.

Now, we first find an odd ``fat cycle'' $C'=v_1'v_2'\dots v_{s-1}'v_s'v_1'$ of length $s$ in $G_\lambda$, where $v_i'\in V_i'$.
Then, from the above observation, we can enlarge this fat cycle $C'$ to an odd cycle of length from $s$ to
\[
(s-1)(m-5\epsilon m/\beta).
\]
Note that (\ref{eta-ep}) and (\ref{alp}), it follows that
\begin{align*}
(s-1)(m-5\epsilon m/\beta)&\ge \left(\frac{(1+0.1\eta)t}{4(1+\eta)}-1\right)\left(1-\frac{5\epsilon} {2p/3-8\epsilon}\right)(1-\epsilon)^2\frac{N}{t}
\\&>\frac{t}{4(1+\eta)}(1-5\sqrt{\epsilon})(1-2\epsilon)\frac{(8+\alpha)n}{t}
\\&>(2+\eta)n.
\end{align*}
Therefore, by noting that $s\le t$ which is much smaller than $n_\lambda$, we can find an odd cycle of length $2n_\lambda+1$ as desired. This completes the proof of the three color case.

\medskip
For the two color case, we use a result by Nikiforov and Schelp \cite[Theorem 3]{niki} instead of Lemma \ref{Luczak-1999}.

\begin{lemma}(Nikiforov and Schelp \cite{niki})
If $t$ is large and $H$ is a graph of order $2t-1$, with minimum degree $\delta(H)\ge(2-10^{-6})t$, then for every 2-coloring of $E(H)$ one of the colors contains cycles $C_s$ for all $s\in[3,t]$.
\end{lemma}

Since a graph $H$ with at least $(1-\epsilon){t\choose 2}$ edges contains a subgraph of order at least $(1-\sqrt{\epsilon})t$ with minimum degree at least $(1-3\sqrt{\epsilon})t$.
Then a similar proof for the two color case follows.
This completes the proof of Theorem \ref{odd-c-3}.
\hfill $\Box$

\medskip
In the following, we focus on the even cycle Ramsey number of random graph. We will give a proof since there are still many details to handle although the main idea is similar to that of Theorem \ref{odd-c-3}.
Let us first recall the following result by Figaj and {\L}uczak \cite[Lemma 8]{Figaj--2007}.

\begin{lemma}(Figaj and {\L}uczak \cite{Figaj--2007})\label{fl-2007}
 For every $0<\alpha_3\le\alpha_2\le\alpha_1=1$ and $0<\epsilon<10^{-21}\alpha^7_3$,
there exists $t_0$, such that for each graph $H$ with $T=(2+\alpha_2+\alpha_3+18\epsilon^{1/7})t\ge t_0$
 vertices and at least $(1-\epsilon){T\choose 2}$ edges, the following holds. For every 3-coloring of
the edges of $H$ there exists a color $i$, $i = 1, 2, 3$, such that some component of the subgraph
induced in $H$ by the edges of the $i$th color contains a matching saturating at least $(2\alpha_i+0.1\epsilon^{1/7})t$
vertices.
\end{lemma}

In order to prove Theorem \ref{c-3}, it is equivalent to prove the following theorem.
\begin{theorem}\label{c-3-1}
Let $0<p\le1$ be fixed. For $\alpha_1\ge\alpha_2\ge\alpha_3>0$, we have a.a.s. that
\begin{align*}
\mathcal{G}((2\alpha_1+\alpha_2+\alpha_3+o(1))n,p) &\to (C_{2\lfloor\alpha_1n\rfloor},C_{2\lfloor\alpha_2n\rfloor},C_{2\lfloor\alpha_3n\rfloor}),\;\;\mbox{and}\\
\mathcal{G}((2\alpha_1+\alpha_2+o(1))n,p) &\to (C_{2\lfloor\alpha_1n\rfloor},C_{2\lfloor\alpha_2n\rfloor}).
\end{align*}
\end{theorem}

\noindent
{\em 3.2 \; Proof of Theorem \ref{c-3-1}}

\bigskip

We first prove the three color case.
Let $0<p\le1$ be fixed,
\begin{align}\label{ep-t0}
0<\epsilon\le\min\left\{10^{-21}\alpha^7_3,{p^2}/{4}\right\},
\end{align}
and $t_0$ be defined as in Lemma \ref{fl-2007}.
Furthermore, let $\mu=2\alpha_1+\alpha_2+\alpha_3+\delta$,
\begin{align}\label{alp-n}
\delta=100\epsilon^{1/7}\alpha_1/\alpha_3,\;\;\text{and}\;\;N=\mu n.
\end{align}
We will show that a.a.s. for every 3-edge coloring of $G=\mathcal{G}(N,p)$ with vertex set $V$, there is a monochromatic even cycle of length $2\lfloor\alpha_\theta n\rfloor$ for some $\theta=1,2,3$.

Consider a 3-edge coloring $(G_1, G_2, G_3)$ of $G$.
By Lemma \ref{reg}, there exists $T_0=T_0(\epsilon,t_0)$ such that there is a partition $\{V_0,V_1,\dots,V_T\}$ of $V(G)$ satisfying $t_0\le T\le T_0$ and
\medskip

(1) $|V_0|<\epsilon N, |V_1|=|V_2|=\dots=|V_T|$;

(2) all but at most $\epsilon T^2$ pairs $(V_i,V_j)$, $i,j\in [T]$, are $\epsilon$-regular for $G_1$, $G_2$ and $G_3$.

\medskip

Similar to Claim \ref{claim-1}, we suppose that $d_{G}(V_i,V_j)\geq p/2$ for every  $1\leq i \leq j \leq T$.
Let $H$ be the auxiliary graph on vertex set $[T]$,
where $ij$ is an edge of $H$ if and only if  $(V_i,V_j)$ is $\epsilon$-regular.
Since the partition $V_0,V_1,\dots,V_T$ is $\epsilon$-regular, the  number of edges in $H$
is at least $(1-\epsilon){T\choose 2}.$
We color an edge $ij$ in $H$ red if $d_{G_1}(V_i,V_j)\geq p/6$, blue if $d_{G_2}(V_i,V_j)\geq p/6$ and green if $d_{G_3}(V_i,V_j)\geq p/6$.

Let $\alpha_\theta'=\alpha_\theta/\alpha_1$ for $\theta=1,2,3$, and let $t$ be the maximum integer such that
\begin{align*}
T\ge(2+\alpha_2'+\alpha_3'+18\epsilon^{1/7}) t,
\end{align*}
which implies that
\begin{align}\label{T-t}
t\ge(2+\alpha_2'+\alpha_3'+18\epsilon^{1/7})^{-1}T-1.
\end{align}

From Lemma \ref{fl-2007},
there exists a color $\theta$, $\theta = 1, 2, 3$, such that some component of the subgraph
induced in $H$ by the edges of the $\theta$th color contains a matching saturating at least $(2\alpha'_\theta+0.1\epsilon)t$
vertices. Thus,
we can obtain a monochromatic path $P$ of even length $s$ on at
least $(2\alpha'_\theta+0.1\epsilon)t$ vertices.
Without loss of generality, suppose that the path $P$ has vertex set $[s]$, i.e. for $i=1,2,\dots,s-1$,

\medskip
(a) $(V_i,V_{i+1})$ is $\epsilon$-regular;

\smallskip
(b) $d_{G_\theta}(V_i,V_{i+1})\geq p/6$.

\medskip
In the following, we aim to show that $G_1$ contains an even cycle of length $2\lfloor\alpha_\theta n\rfloor$.
To this end, for $1\le i\le s$, we divide  each set $V_{i}$ into two subsets $U_i$, $W_i$ of sizes as equal as possible,
that is, $\lfloor\frac{(1-\epsilon)N}{2T}\rfloor\le|U_{i}|,|W_{i}|\le \lfloor\frac{(1-\epsilon)N}{2T}\rfloor+1.$ We delete the ceiling and floor signals as it dose not affect the final result. So we admit that $|U_{i}|=|W_{i}|=\frac{(1-\epsilon)N}{2T}$.
For $1\le j\le s/2$, by Fact \ref{hdt}, we have $(U_{2j-1},U_{2j})$ is $2\epsilon$-regular with density at least $(p/6-\epsilon)$.
Now, Fact \ref{supreg} implies that we can take $U_{2j-1}'\subseteq U_{2j-1}$ and $U_{2j}'\subseteq U_{2j}$
with $|U_{2j-1}'|\ge (1-2\epsilon)|U_{2j-1}|$ and $|U_{2j}'|\ge (1-2\epsilon)|U_{2j}|$ such that
$(U'_{2j-1},U'_{2j})$ is $(4\epsilon,p/6-5\epsilon)$-super-regular for $G_\theta$.
Similarly, for $2\le j\le s/2$, we can take $W_{2j-1}'\subseteq W_{2j-1}$ and $W_{2(j-1)}'\subseteq W_{2(j-1)}$
with $|W_{2j-1}'|\ge (1-2\epsilon)|W_{2j-1}|$ and $|W_{2(j-1)}'|\ge (1-2\epsilon)|W_{2(j-1)}|$ such that $(W'_{2j-1},W'_{2(j-1)})$ is $(4\epsilon,p/6-5\epsilon)$-super-regular for $G_\theta$.
Let
\begin{align*}
m=(1-2\epsilon)(1-\epsilon)\frac{N}{2T},\;\;\epsilon'=4\epsilon,\;\;\text{and}\;\;\beta=2p/3-20\epsilon.
\end{align*}
Therefore, by Lemma \ref{ben-sk}, for every $\ell$ with $1\le \ell \le m-5\epsilon m/\beta$,
and for every pair of vertices $u_{2j-1}'\in U_{2j-1}'$ and $u_{2j}'\in U_{2j}'$, $G_\theta$ contains a path of length $2\ell+1$ connecting $u_{2j-1}'$ and $u_{2j}'$ for $1\le j\le s/2$.
Similarly, for every pair of vertices $w_{2j-1}'\in W_{2j-1}'$ and $w_{2(j-1)}'\in W_{2(j-1)}'$, $G_\theta$ contains a path of length $2\ell+1$ connecting $w_{2j-1}'$ and $w_{2j}'$ for $2\le j\le s/2$.

Now, we first find a ``fat path'' $P'=u_1'u_2'u_3'u_4'\dots u_{s-1}'u_s'w_{s-1}'w_{s-2}'\dots w_3'w_2'$ of length $2(s-1)$ in $G_\theta$, where $u_i'\in U_i'$ and $w_i'\in W_i'$.
From the above observation, we can enlarge this fat path $P'$ to an odd path of length from $2s-3$ to
\[
(2s-3)(m-5\epsilon m/\beta).
\]
In order to get an even cycle of length $\lfloor2\alpha_\theta n\rfloor$ as desired, let us fix such a path $Q$ such that both the ``fat edges'' $u_1'u_2'$ and $w_3'w_2'$ contain at least $2\epsilon N$ vertices.
Then
\[
|V(Q)\cap V_1|\ge\epsilon N\;\; \text{and} \;\;|V(Q)\cap V_2|\ge\epsilon N.
\]
Since $(V_1,V_2)$ is $\epsilon$-regular, we can find an edge between $V(Q)\cap V_1$ and $V(Q)\cap V_2$.
Consequently, we can find even cycles of length from $2s-2+2\epsilon N$ to
\[
n'=(2s-3)(m-5\epsilon m/\beta)+1-2\epsilon N.
\]
Note that (\ref{ep-t0}) and (\ref{alp-n}), it follows that
\begin{align*}
n'&\ge (2s-3)\left(1-5\epsilon /\beta\right) (1-2\epsilon)(1-\epsilon)\frac{N}{2T}-2\epsilon N
\\&>2(2\alpha'_\theta+0.05\epsilon^{1/7})(2+\alpha_2'+\alpha_3'+18\epsilon^{1/7})^{-1}T
\cdot(1-5\sqrt{\epsilon})(1-3\epsilon)\frac{\mu n}{2T}-2\epsilon N
\\&>(2\alpha'_\theta+0.05\epsilon^{1/7})(2+\alpha_2'+\alpha_3'+18\epsilon^{1/7})^{-1}(1-8\sqrt{\epsilon})\mu n-2\epsilon N
\\&>(2\alpha_\theta+0.05\epsilon^{1/7})n.
\end{align*}
Therefore, we can find an even cycle of length $2\lfloor\alpha_\theta n\rfloor$ as desired.
This completes the proof of the three color case.

\medskip
For the two color case, we use a result by Letzer \cite[Theorem 2]{letzer-2016} instead of Lemma \ref{fl-2007}.

\begin{lemma}(Letzer \cite{letzer-2016})\label{letz}
Let $0<\epsilon<1/4$, $t_1\ge t_2$, and let $H$ be a graph of order $T\ge t_1+\lfloor (t_2+1)/2\rfloor+150\sqrt{\epsilon}t_1$
 with at least $(1-\epsilon){T\choose 2}$ edges. Then $H\to (P_{t_1},P_{t_2})$.
\end{lemma}

A similar proof for the two color case follows by using Lemma \ref{letz}.
This completes the proof of Theorem \ref{c-3}.
\hfill $\Box$

\section{Concluding remarks}

Let us point out that the assertions we obtained are valid when the probability $p$ is fixed.
It would be interesting that whether the assertions hold also as that for path when $pn\to\infty$,
or somewhat weakly that whether the assertions hold when $p\to0$.
As an end, we propose the following problems.

\begin{problem}
Prove or disprove as $pn\to\infty$, we have a.a.s. that
\begin{align*}
\mathcal{G}((8+o(1))n,p) \to (C_{2n+1})_3,\;\;&\text{and}\;\;
\mathcal{G}((4+o(1))n,p) \to (C_{2n+1})_2;\\
\mathcal{G}((4+o(1))n,p) \to (C_{2n})_3,\;\;&\text{and}\;\;
\mathcal{G}((3+o(1))n,p) \to (C_{2n})_2.
\end{align*}
\end{problem}

A weaker type is as follows.

\begin{problem}
Prove or disprove as $p\to0$, we have a.a.s. that
\begin{align*}
\mathcal{G}((8+o(1))n,p) \to (C_{2n+1})_3,\;\;&\text{and}\;\;
\mathcal{G}((4+o(1))n,p) \to (C_{2n+1})_2;\\
\mathcal{G}((4+o(1))n,p) \to (C_{2n})_3,\;\;&\text{and}\;\;
\mathcal{G}((3+o(1))n,p) \to (C_{2n})_2.
\end{align*}
\end{problem}

\end{spacing}


\begin{thebibliography}{99}

\bibitem{Alon-1992}
N. Alon and J. Spencer, The Probabilistic Method, Wiley--Interscience, New York, 1992.



\bibitem{Benevides-2009}
F. S. Benevides and J. Skokan, The 3-colored Ramsey number of even cycles, J. Combin. Theory
Ser. B 99 (2009), 690--708.


\bibitem{B-e}
J.A. Bondy, P. Erd\H{o}s, Ramsey numbers for cycles in graphs, J. Combin. Theory Ser. B 14 (1973), 46--54.

\bibitem{Chernoff-1952}
H. Chernoff, A measure of the asymptotic efficiency for tests of a hypothesis based on the sum of observations,
Ann. Math. Statis 23 (1952), 493--509.

\bibitem{cfs}
D. Conlon, J. Fox, and B. Sudakov, Recent developments in graph
Ramsey theory. Surveys in Combinatorics 2015. Edited by Artur Czumaj,
Agelos Georgakopoulos, Daniel Kr¨¢l, Vadim Lozin, Oleg Pikhurko.
Cambridge University Press, pp 49--118.

\bibitem{Davies-2017}
E. Davies, M. Jenssen, and B. Roberts, Multicolour Ramsey numbers of paths and even cycles,
European J. Combin. 63 (2017), 124--133.

\bibitem{dell}
D. Dellamonica Jr., Y. Kohayakawa, M. Marciniszyn, and A. Steger, On the resilience
of long cycles in random graphs, Electron. J. Combin. 15 (2008), R32.


\bibitem{dudek-2018}
A. Dudek and P. Pra{\l}at, On some multicolour Ramsey properties of random graphs, SIAM J. Discrete Math. 31 (2017), 2079--2092.



\bibitem{er}
P. Erd\H{o}s, R.J. Faudree, C.C. Rousseau, R.H. Schelp, Generalized Ramsey theory for multiple colors,
J. Combin. Theory Ser. B 20 (1976), 250--264.

\bibitem{Faudree-1974}
R. Faudree and R.Schelp, All Ramsey numbers for cycles in graphs, Discrete Math. 8 (1974), 313--329.

\bibitem{Figaj--2007}
A. Figaj and T. {\L}uczak, The Ramsey number for a triple of long even cycles,  J. Combin.
Theory Ser. B 97 (2007), 584--596.

\bibitem{f-l}
A. Figaj and T. {\L}uczak, The Ramsey number for a triple of long cycles, Combinatorica, https://doi.org/10.1007/s00493-016-2433-y.

 \bibitem{graham-roth-spen}
R. Graham, B. Rothschild and J. Spencer, Ramsey Theory, Wiley, New
York, 1980.


\bibitem{Gyarfas-2007}
A. Gy\'{a}rf\'{a}s, M. Ruszink\'{o}, G. N. S\'{a}rk\"{o}zy, and E. Szemer\'{e}di, Three-color Ramsey numbers for
paths, Combinatorica 27 (2007),  35--69.

\bibitem{Gyarfas-M-2007}
A. Gy\'{a}rf\'{a}s, M. Ruszink\'{o}, G. N. S\'{a}rk\"{o}zy, and E. Szemer\'{e}di, Tripartite Ramsey numbers
for paths, J. Graph Theory 55 (2007), 164--174.



\bibitem{Jenssen-Skokan-2016}
M. Jenssen and J. Skokan, Exact Ramsey numbers of odd cycles via nonlinear optimisation,
arXiv:1608.05705 (2016), preprint.

\bibitem{k-r}
G. Karolyi and V. Rosta, Generalized and geometric Ramsey numbers for cycles,
Theoret. Comp. Sci. 263 (2001), 87--98.

\bibitem{Ko-r}
Y. Kohayakawa, V. R\"{o}dl, Szemer\'{e}di's regularity lemma and quasi-randomness, in: Recent Advances in Algorithms
and Combinatorics, in: CMS Books Math./Ouvrages Math. SMC, vol. 11, Springer, New York, 2003, 289--351.

\bibitem{Kohayakawa-2005}
Y. Kohayakawa, M. Simonovits, and J. Skokan, The 3-colored Ramsey number of odd cycles,
Electron. Notes Discrete. Math. 19 (2005), 397--402.

\bibitem{koha97}
Y. Kohayakawa, Szemer\'edi's Regularity Lemma for sparse graphs, In Foundations of Computational Mathematics: Rio de Janeiro, 1997, Springer, (1997), 216--230.





\bibitem{Kom-Simonovits-2005}
J. Koml\'{o}s and M. Simonovits, Szemer\'{e}di's Regularity Lemma and its applications in graph
theory, Combinatorics, Paul Erd\H{o}s is Eighty, vol. 2, Bolyai Society Mathematical Studies, 1996,
295--352.

\bibitem{Knierim-2018}
C. Knierim and P. Su, Improved bounds on the multicolor Ramsey numbers of paths and even cycles,
arXiv: 1801.04128v1.

\bibitem{letzer-2016}
S. Letzter: Path Ramsey number for random graphs, Combin.  Probab. Comput. 25 (2016), 612--622.

\bibitem{Luczak-1999}
T. {\L}uczak, $R(C_n,C_n,C_n)\leq (4+o(1))n$, J. Combin. Theory Ser. B 75 (1999), 174--187.

\bibitem{niki}
V. Nikiforov and R. Schelp, Cycles and Stability, J. Combin. Theory Ser. B 98 (2008), 69--84.


\bibitem{Rosta-1973}
V. Rosta, On a Ramsey-type problem of J. A. Bondy and P. Erd\H{o}s. I, II, J. Combin. Theory
Ser. B 15 (1973), 105--120.


 \bibitem{roth53}
K. Roth, On certain sets of integers, J. London Math. Soc. 28 (1953),  104-109.

\bibitem{roth54}
K. Roth, On certain sets of integers, II, J. London Math. Soc. 29 (1954),  20-26.


\bibitem{Gsark-2016}
G. S\'{a}rk\"{o}zy, On the multi-coloured Ramsey numbers of paths and even cycles, Electron. J. Combin.
23 (2016), P3.53.

\bibitem{sun}
Y. Sun, Y. Yang, F. Xu, and B. Li, New lower bounds on the multicolor Ramsey
numbers $R_r(C_{2m})$. Graphs and Combinatorics 22 (2006), 283--288.

\bibitem{sze75}
E. Szemer\'edi, On sets of integers containing no $k$ elements in arithmetic progression, Acta Arithmetica 27 (1975),  199-245.


\bibitem{Sze-1978}
E. Szemer\'{e}di,   Regular partitions of graphs. Probl\`{e}mes Combinatoires et Th\'{e}orie des
Graphes: Colloq. Internat. CNRS, Univ. Orsay, Orsay, 1976, Vol. 260 of Colloq. Internat. CNRS,
CNRS, Paris,  399--401.




\end{thebibliography}
\end{document}